\documentstyle[amssymb, 10pt]{article}

\begin{document}
\font\frak=eufm10 scaled\magstep1
\font\fak=eufm10 scaled\magstep2
\font\fk=eufm10 scaled\magstep3
\font\scriptfrak=eufm10
\font\tenfrak=eufm10
\font\msb=eusb10

\newtheorem{theorem}{Theorem}
\newtheorem{corollary}{Corollary}
\newtheorem{proposition}{Proposition}
\newtheorem{definition}{Definition}
\newtheorem{lemma}{Lemma}
\font\frak=eufm10 scaled\magstep1
\newenvironment{pf}{{\noindent{\it Proof. }}}{\ $\Box$\medskip}


\mathchardef\za="710B  
\mathchardef\zb="710C  
\mathchardef\zg="710D  
\mathchardef\zd="710E  
\mathchardef\zve="710F 
\mathchardef\zz="7110  
\mathchardef\zh="7111  
\mathchardef\zvy="7112 
\mathchardef\zi="7113  
\mathchardef\zk="7114  
\mathchardef\zl="7115  
\mathchardef\zm="7116  
\mathchardef\zn="7117  
\mathchardef\zx="7118  
\mathchardef\zp="7119  
\mathchardef\zr="711A  
\mathchardef\zs="711B  
\mathchardef\zt="711C  
\mathchardef\zu="711D  
\mathchardef\zvf="711E 
\mathchardef\zq="711F  
\mathchardef\zc="7120  
\mathchardef\zw="7121  
\mathchardef\ze="7122  
\mathchardef\zy="7123  
\mathchardef\zf="7124  
\mathchardef\zvr="7125 
\mathchardef\zvs="7126 
\mathchardef\zf="7127  
\mathchardef\zG="7000  
\mathchardef\zD="7001  
\mathchardef\zY="7002  
\mathchardef\zL="7003  
\mathchardef\zX="7004  
\mathchardef\zP="7005  
\mathchardef\zS="7006  
\mathchardef\zU="7007  
\mathchardef\zF="7008  
\mathchardef\zW="700A  

\newcommand{\be}{\begin{equation}}
\newcommand{\ee}{\end{equation}}
\newcommand{\ra}{\rightarrow}
\newcommand{\lra}{\longrightarrow}
\newcommand{\bea}{\begin{eqnarray}}
\newcommand{\eea}{\end{eqnarray}}
\newcommand{\beas}{\begin{eqnarray*}}
\newcommand{\eeas}{\end{eqnarray*}}
\newcommand{\Z}{{\Bbb Z}}
\newcommand{\R}{{\Bbb R}}
\newcommand{\C}{{\Bbb C}}
\newcommand{\1}{{\bold 1}}
\newcommand{\SL}{SL(2,\C)}
\newcommand{\Sl}{sl(2,\C)}
\newcommand{\SU}{SU(2)}
\newcommand{\su}{su(2)}
\newcommand{\G}{{\goth g}}
\newcommand{\D}{{\rm d}}
\newcommand{\de}{\,{\stackrel{\rm def}{=}}\,}
\newcommand{\we}{\wedge}
\newcommand{\nn}{\nonumber}
\newcommand{\ot}{\otimes}
\newcommand{\s}{{\textstyle *}}
\newcommand{\ts}{T^\s}
\newcommand{\da}{\dagger}
\newcommand{\pa}{\partial}
\newcommand{\ti}{\times}
\newcommand{\A}{{\cal A}}
\newcommand{\Li}{{\cal L}}
\newcommand{\ka}{{\Bbb K}}
\newcommand{\find}{\mid}
\newcommand{\ad}{{\rm ad}}
\newcommand{\rS}{]^{SN}}
\newcommand{\rb}{\}_P}
\newcommand{\p}{{\sf P}}
\newcommand{\h}{{\sf H}}
\newcommand{\X}{{\cal X}}
\newcommand{\I}{\,{\rm i}\,}
\newcommand{\rB}{]_P}
\newcommand{\Ll}{\Li}
\def\lna{\lbrack\! \lbrack}
\def\rna{\rbrack\! \rbrack}
\def\rnaf{\rbrack\! \rbrack^\zF}
\def\rnah{\rbrack\! \rbrack\,\hat{}}
\def\lbo{{\lbrack\!\!\lbrack}}
\def\rbo{{\rbrack\!\!\rbrack}}
\newcommand{\Llf}{\Li^\zF}
\newcommand{\Llh}{{\pounds}}
\newcommand{\Lih}{\hat\Li}
\def\tU{\tilde U}
\def\hU{\hat U}
\def\Df{\D^\zF}
\def\Dh{\hat\D}
\def\vm{{\vert\zm\vert}}
\def\Di{{\cal D\,}^\za}
\def\De{{\cal DE\,}^\za}
\def\b{\bar}
\def\la{\langle}
\def\ran{\rangle}
\def\U{{\bf U}}

\title{The graded Jacobi algebras and (co)homology}

\author{
Janusz Grabowski\thanks{Supported by KBN, grant No. 2 P03A 031 17.}\\
Mathematical Institute, Polish Academy of Sciences\\
ul. \'Sniadeckich 8, P. O. Box 137, 00-950 Warszawa, Poland\\
{\it e-mail:} jagrab@impan.gov.pl
\and
Giuseppe Marmo\thanks{Supported by PRIN SINTESI}\\
Dipartimento di Scienze Fisiche,
Universit\`a Federico II di Napoli\\
and\\
INFN, Sezione di Napoli\\
Complesso Universitario di Monte Sant'Angelo\\
Via Cintia, 80126 Napoli, Italy\\
{\it e-mail:} marmo@na.infn.it}
\maketitle
\begin{abstract}\noindent Jacobi algebroids (i.e. `Jacobi versions'
of Lie algebroids)
are studied in the context of graded Jacobi brackets on graded
commutative algebras. This unifies varios concepts of graded Lie
structures in geometry and physics. A method of describing such
structures by classical Lie algebroids via certain gauging (in the
spirit of E.Witten's gauging of exterior derivative) is developed.
One constructs a corresponding Cartan differential calculus
(graded commutative one) in a natural manner. This, in turn, gives
canonical generating operators for triangular Jacobi algebroids.
One gets, in particular, the Lichnerowicz-Jacobi homology
operators associated with classical Jacobi structures.
Courant-Jacobi brackets are obtained in a similar way and use to
define an abstract notion of a Courant-Jacobi algebroid and
Dirac-Jacobi structure.
\end{abstract}

\section{Introduction}
In this paper we propose a unification of various concepts of graded
brackets  one  meets  in geometry  and
physics and a method of `gauging' which allows to pass from the world of
derivations (i.e. tangent bundles, vector fields,  exterior derivatives,
Lie algebroids) to the world of first-order differential operators in
the  spirit  of E.~Witten's \cite{Wi} gauging of  exterior  derivative.
This algebra and geometry is noncommutative in the sense that bosonic
and fermionic parts are incorporated in a unique scheme.
We concentrate on purely mathematical aspects to keep the  size  of  the
paper readable but we hope that possible applications to Batalin-Vilkovisky
formalism, BRST-method, integrability and Dirac structures,  etc., will
be found.

\medskip
For  a  vector  bundle   $E$   over   the   base   manifold   $M$,   let
$A(E)=\oplus_{k\in\Z}A^k(E)$ be the exterior algebra of multisections of
$E$.  This  is  a  basic  geometric  model  for  a  graded   associative
commutative algebra with unity. We will
refer to elements of $\zW^k(E)=A^k(E^*)$ as to {\it k-forms} on $E$.
Here, we identify $A^0(E)=\zW^0(E)$  with  the  algebra  $C^\infty(M)$  of
smooth functions on the base and $A^k(E)=\{ 0\}$ for $k<0$.
Denote by $\vert X\vert$ the Grassmann degree of the multisection $X\in
A(E)$.

As it has been observed in \cite{KS}, a {\it Lie algebroid} structure on
$E$ (for the traditional definition and main properties we recommend the
survey article \cite{Ma}) can be identified  with  a  {\it  Gerstenhaber
algebra} structure (in the terminology of \cite{KS}) on $A(E)$ which  is
just a {\it graded Poisson bracket} on $A(E)$ of degree -1. Recall  that
a graded Poisson bracket of degree  $k$  on  a  $\Z$-graded  associative
commutative algebra $\A=\oplus_{i\in\Z}\A^i$ is a graded bilinear map
\be
\{\cdot,\cdot\}:\A\ti\A\ra\A
\ee
of degree $k$ (i.e. $\vert\{  a,b\}\vert=\vert  a\vert+\vert  b\vert+k$)
such that
\begin{enumerate}
\item $\{ a,b\}=-(-1)^{(a+k)(b+k)}\{ b,a\}$ (graded anticommutativity),
\item $\{ a,bc\}=\{  a,b\}  c+(-1)^{(a+k)b}b\{  a,c\}$  (graded  Leibniz
rule),
\item $\{\{ a,b\},c\}=\{ a,\{  b,c\}\}-(-1)^{(a+k)(b+k)}\{ b,\{  a,c\}\}$
(graded Jacobi identity).
\end{enumerate}
Here we use  the convention that we write just $a$ for $\vert
a\vert$.

\medskip
It is obvious  that  this  notion  extends  naturally  to  more  general
gradings in the algebra. The Leibniz rule tells that the Poisson bracket
is identified with a skew biderivation $\zL$ on $\A$ (`bivector
field'), $\{ a,b\}=\zL(a,b)$, which, due to the Jacobi identity,  has  a
special property (`the Schouten bracket $\lna\zL,\zL\rna$ vanishes').
A precise meaning for the notions in the quotation marks has been  given
in \cite{Kr2}. With respect to  the  standard  terminology,  the  graded
derivation $\zL_a=\{ a,\cdot\}$ associated with $a\in\A$ is  called  the
corresponding {\it hamiltonian vector field} and the map $a\mapsto\zL_a$
is a  homomorphism  from  $(\A,\{\cdot,\cdot\})$  into  the  graded  Lie
algebra $(Der(\A),\lbo\cdot,\cdot\rbo)$ of graded  derivations  of  $\A$
with the graded commutator $\lbo\cdot,\cdot\rbo$.

For  a  graded  commutative  algebra   with   unity   $\1$,   a   natural
generalization of  a  graded  Poisson  bracket  is  {\it  graded  Jacobi
bracket}. The only difference is that we replace the Leibniz rule by
\be\label{first}
\{ a,bc\}=\{  a,b\}  c+(-1)^{(a+k)b}b\{  a,c\}-\{ a,\1\} bc
\ee
which just means that  $\{  a,\cdot\}$  is  a  first-order  differential
operator on $\A$ (for the differential calculus  on  graded  commutative
algebras  we  refer  to   \cite{ViA,VK}).   This   goes   back   to   the
well-known observation by Kirillov \cite{Ki} that in the case  of
$\A=C^\infty(M)$ every local Lie bracket is of first order (an algebraic
version of  this fact in ungraded case has been proved in \cite{Gr}).
A graded associative commutative algebra with a graded Jacobi  structure
we will call a {\it graded Jacobi algebra}.

\bigskip\noindent
{\bf Definition.} A {\it  graded  Jacobi algebra} is a graded, say
$\Z^n$-graded, associative  commutative  algebra  $\A=\oplus_{i\in\Z^n}$
with  unit  $\1$  equipped  with   a   {\it   graded   Jacobi   bracket}
$\{\cdot,\cdot\}$  of degree $k\in\Z^n$, i.e. a graded bilinear map
\be
\{\cdot,\cdot\}:\A\ti\A\ra\A
\ee
of degree $k$ (i.e. $\vert\{  a,b\}\vert=\vert  a\vert+\vert  b\vert+k$)
such that
\begin{enumerate}
\item $\{ a,b\}=-(-1)^{\langle a+k,b+k\rangle}\{ b,a\}$ (graded
anticommutativity),
\item $\{ a,bc\}=\{   a,b\}   c+(-1)^{\langle  a+k,b\rangle}b\{   a,c\}
-\{a,\1\}bc$  (graded generalized Leibniz rule),
\item $\{\{ a,b\},c\}=\{ a,\{   b,c\}\}-(-1)^{\langle  a+k,b+k\rangle}\{
b,\{  a,c\}\}$ (graded Jacobi identity),
\end{enumerate}
where $\langle\cdot,\cdot\rangle$ is the standard pairing in $\Z^n$.

\bigskip\noindent
The generalized Leibniz rule tells that the bracket is a  bidifferential
operator on $\A$ of first order.  In  the  non-graded  case,  under  the
assumption that there are no  non-trivial  nilpotent  elements  in  $\A$,
every Lie bracket given by bidifferential operator  is  known  \cite{Gr}
to be of first  order  (it  is  a  generalization  of  this  result  for
$\A=C^\infty(M)$ \cite{Ki}).
In the classical case of the algebra $C^\infty(M)$, every  skew-symmetric
first-order bidifferential operator  $J$  splits  into  $J=\zL+I\we\zG$,
where $\zL$ is a bivector field, $\zG$ is a  vector  field  and  $I$  is
identity, so that the corresponding bracket of functions reads
\be\label{jac}
\{ f,g\}_J=\zL(f,g)+f\zG(g)-g\zG(f).
\ee
The Jacobi identity for this bracket is usually written in terms of  the
Schouten-Nijenhuis bracket by
\bea
\lna\zG,\zL\rna&=&0,\\
\lna\zL,\zL\rna&=&-2\zG\we\zL.
\eea
Hence, every Jacobi bracket on $C^\infty(M)$ can be identified with  the
pair $J=(\zL,\zG)$ satisfying the above conditions,  i.e.  with  a  {\it
Jacobi structure} on $M$ (cf. \cite{Li}).
Note that we use the version of  the  Schouten-Nijenhuis  bracket  which
gives  a  graded  Lie  algebra  structure  on  multivector  fields  (cf.
\cite{Mi}) and which differs from the classical one by signs.
On the other hand, this Jacobi identity  can be written in terms of  the
algebraic   Richardson-Nijenhuis   bracket (cf. \cite{NR}) $\lna
J,J\rna^{RN}=0$ of skew-multiplicative maps on $C^\infty(M)$  which,  as
it has been observed in \cite{GM}, is a deformation of a Schouten-type
bracket,  when    reduced    to    first-order   {\it   polydifferential
operators}  on  $C^\infty(M)$,  i.e.  skew-symmetric   multidifferential
operators.  However,  this  bracket  on   first-order   polydifferential
operators is not the Schouten  bracket  for   the   Lie   algebroid   of
linear scalar first-order differential operators but a Jacobi version of
it. This means that here we have  a difference similar to the difference
between Poisson and  Jacobi  brackets  in  the  classical  case,  if  we
understand the Schouten bracket as being a graded Poisson bracket
according to \cite{Kr2}.

An  analogous  construction  for  a  general  Lie  algebroid  has   been
introduced in \cite{IM} under  the  name  of  a  {\it  generalized  Lie
algebroid}. We have recognized this structure as  being  an  odd  Jacobi
structure in \cite{GM} (in fact, a graded Jacobi bracket of degree -1 in
the terminology of this paper), where the name {\it Jacobi algebroid}
has  been used. These structures are closely related to the  concept  of
Lie-like brackets on affine bundles \cite{GGU,MMS} as it has been
mentioned in \cite{GGU}.

The {generalized Lie algebroids} in the terminology of \cite{IM}  or
the {Jacobi algebroids} in the terminology of \cite{GM} are
associated with the
pairs combined with a Lie algebroid bracket on a vector bundle $E$ over
$M$ and a 1-cocycle $\zF\in\zW^1(E)$, $\D\zF=0$, relative to the Lie
algebroid exterior derivative $\D$. The Schouten-Jacobi bracket on the
graded algebra $A(E)$ of multisections of $E$ is given by
\be\label{jb}
\lna X,Y\rnaf=\lna X,Y\rna+xX\we i_\zF Y-(-1)^xyi_\zF X\we Y       ,
\ee
where we use the convention  that  $x=\vert  X\vert-1$  is  the  shifted
degree of $X$ in the graded algebra $A(E)$ and $\lna\cdot,\cdot\rna$  is
the  Schouten  bracket  of  the  corresponding  Lie  algebroid.
Note that $i_\zF X=(-1)^x\lna X,\1\rnaf$ and (\ref{first}) is satisfied:
\be\label{first1}
\lna X,Y\we Z\rnaf=\lna X,Y\rnaf\we Z+(-1)^{x(y+1)}Y\we\lna X,Z\rnaf -
\lna X,\1\rnaf\we Y\we Z.
\ee
The Jacobi algebroid bracket (\ref{jb}) can be written in the form
\bea\nn
\lna X,Y\rnaf&=&\lna X,Y\rna+\vert X\vert X\we i_\zF  Y-(-1)^x
\vert Y\vert i_\zF X\we Y
+(-1)^xi_\zF(X\we Y)=\\ \label{gm}
&&\lna X,Y\rna+(Deg\we i_\zF)(X,Y)-X\we i_\zF Y+(-1)^xi_\zF X\we Y,
\eea
where $Deg(X)=\vert X\vert X$ is a derivation of $A(E)$ and $Deg\we
i_\zF$ is an
appropriate wedge product of graded derivations. This  shows  that  the
above bracket is essentially like (\ref{jac}) with $\zL=\lna\cdot, \cdot
\rna+Deg\we i_\zF$ and $\zG=-i_\zF$.

Note that the bracket $\lna\cdot,\cdot\rnaf$ is completely determined by
its values on the  Lie  subalgebra  $C^\infty(M)\oplus\textrm{Sec}(E)$
due to the generalized Leibniz property (\ref{first1}). The Lie
algebra $C^\infty(M)\oplus\textrm{Sec}(E)$ has the grading induced from
$A(E)$ but  we
must stress that it can be viewed as  a  standard  Lie  algebra  with  a
grading and not a  graded  Lie  algebra  (the  antisymmetry  and  Jacobi
identity is standard not graded). We show in the next section that there
is a {\it linear Jacobi structure} on $E^*$  which  corresponds  to  the
Lie algebra  structure  $(C^\infty(M)\oplus\textrm{Sec}(E),\lna\cdot,
\cdot\rnaf)$. This is a Jacobi analog of the correspondence
\be
\textrm{Lie algebroid structure on }E\ \leftrightarrow\textrm{ linear
Poisson structure on }E^*.
\ee
Since {\it linearity} of tensors on  vector  bundles  is  a   particular
case  of  {\it homogeneity} in the sense of \cite{DLM}  we  will  study
also  relations between homogeneous Poisson and Jacobi structures.

In section 3 we develop a method of inducing Jacobi  algebroid  brackets
by gauging. This is an analog of E.~Witten's  gauging  of  the  exterior
derivative
\be\label{wd}
e^{-f}\D(e^f\zm)=\D\zm+\D f\we\zm
\ee
for  the  Schouten-Nijenhuis  bracket  (the  role  of  the  differential
(\ref{wd}) in studying Jacobi structures has been already observed by
A.~Lichnerowicz \cite{Li}).
We get in this way an appropriate concept of a  Lie  differential  which
gives   generating   operators   of    BV-algebras    associated    with
generalized Jacobi  structures,  and  thus  the  corresponding  homology
operators  for  free.  For  classical  Jacobi  structures  we   end   up
in this way with Lichnerowicz-Jacobi homology and  generating  operators
for associated Lie algebroids (cf. \cite{LLMP,ILMP,Uch,Va1}).

Section 4 is devoted to Courant-Jacobi brackets, i.e.  `Jacobi  versions'
of Courant brackets and Courant algebroids (cf. \cite{Co, LWX, Ro,  Wa,
IM2}). This concept, again, is developed  naturally  by  the  method  of
gauging.

The last section contains pure  algebraic  generalizations  of Jacobi
algebroids. We start with a graded associative commutative algebra  $\A$
and construct the graded Jacobi algebra $\Di(\A)$ of first-order
polydifferential operators on $\A$. The abstract Schouten-Jacobi bracket
on $\Di(\A)$  recognizes  supercanonical  elements  (Jacobi  structures)
which generate graded Jacobi brackets on $\A$, so we  can  consider  the
corresponding cohomology operators as being hamiltonian `vector fields'.
This can be viewed as  a `Jacobi version' of the results of
I.~S.~Krasil'shchik \cite{Kr2} and a
graded algebraic generalization of this kind of structure described in
\cite{GM}. An abstract version of Lichnerowicz-Jacobi cohomology is
defined.

\section{Linear and affine Jacobi structures}

Suppose that we are given a Poisson tensor $\zL$ on a manifold $N$, which
is homogeneous with respect to a vector field $\zD$ (cf. \cite{DLM}), i.e.
$[\zD,\zL]=-\zL$, where $[\cdot,\cdot]$ stands for the Schouten bracket in
the form which gives a graded Lie algebra structure on the graded space
$A(M)=\oplus_iA^i(M)$ of multivector fields on $M$ (this differs by a
sign from the traditional Schouten bracket).

\begin{lemma} The pair $J=(\zL+\zG\we\zD,\zG)$ is a Jacobi structure if
and only if
\be\label{1}
\zG\we [\zD,\zG]\we\zD=[\zG,\zL]\we\zD.
\ee
\end{lemma}
\begin{pf} By direct calculations
\be\label{2}
[\zL+\zG\we\zD,\zL+\zG\we\zD]=[\zL,\zL]-2[\zG,\zL]\we\zD+2\zG\we[\zD,\zL]
+2\zG\we[\zD,\zG]\we\zD.
\ee
Since $[\zL,\zL]=0$ and $[\zD,\zL]=-\zL$, this equals
$-2\zG\we(\zL+\zG\we\zD)=-2\zG\we\zL$ if and only if (\ref{1}) is
satisfied.
\end{pf}

\begin{corollary} If $\zG$ is a homogeneous vector field with
$[\zG,\zL]=0$, then $J=(\zL',\zG)$ with $\zL'=\zL+\zG\we\zD$ is a Jacobi
 structure.
\end{corollary}
\begin{pf}   $\zG$   is   homogeneous   of   degree   $k$   means   that
$[\zD,\zG]=k\zG$ for certain $k$. Thus, $\zG\we[\zD,\zG]=0$ and the
corollary follows by the lemma.
\end{pf}

\begin{corollary} If $f$ is a homogeneous function of degree $k$, then
$J=(\zL',\zG)$ with $\zG=\zL_f$, $\zL'=\zL+\zL_f\we\zD$ and
$\zL_f=i_{df}\zL$ being the hamiltonian  vector  field  associated  with
$f$, is a Jacobi structure. Moreover, if $k\ne 0$, then this structure is
tangent (i.e. $\zL$ and $\zG$ are tangent) to the submanifold determined
by the equation $f=\frac{1}{k}$ (assuming that $\frac{1}{k}$ is a
regular value of $f$).
\end{corollary}
\begin{pf} It is easy to see that the Hamiltonian vector field $\zL_f$ is
homogeneous of degree $(k-1)$. Then, $J$ is a Jacobi structure due to the
previous corollary. The function $f-\frac{1}{k}$ acts by the Jacobi
bracket by
\be
\{f-\frac{1}{k},\cdot\}_J=\zL_f-kf\zL_f+(f-\frac{1}{k})\zL_f=
(f-\frac{1}{k})(1-k)\zL_f
\ee
which vanishes on the described submanifold.
\end{pf}

\medskip\noindent
Note that similar observations have been done in \cite{Pe}.
Many important examples of Jacobi manifolds  are  of the above
form, for instance, contact submanifolds  of exact symplectic  ones  or
spheres in duals of Lie algebras.
A standard example of a homogeneous Poisson tensor is a
linear Poisson tensor on a vector bundle $E$ over a base manifold $M$ with
$\zD$ being the Liouville vector field on $E$. Linearity means that the
Poisson bracket of linear (along fibres) functions on $E$ is again a
linear function. Since every linear function is represented by a section
of the dual bundle $E^*$ by contraction, this gives a Lie bracket on
sections of $E^*$ which is a Lie algebroid bracket giving rise to the
corresponding Schouten bracket on the graded space
$A(E^*)=\oplus_iA^i(E^*)$ of multisections of $E^*$. A generalized
version of this kind of bracket defined in \cite{IM} has been recognized
in \cite{GM} as a (graded) Jacobi bracket of degree -1 on $A(E^*)$. Being
of degree -1, it defines a Lie bracket on a graded subspace
$Aff(E)=C^\infty(M)\oplus\zG(E^*)$ which determines the whole Jacobi
bracket completely, due to being a first order operator. The notation
$Aff(E)$ is justified by the fact that $Aff(E)$ is just the graded space
of affine functions on $E$ with the obvious identification of $C^\infty(M)$
with the algebra of basic functions on $E$. The graded bracket on $Aff(E)$
comes from a homogeneous Jacobi bracket on $E$ determined by the Jacobi
structure $J=(\zL+\zF^v\we\zD,\zF^v)$, where $\zL$ is a linear Poisson
tensor and $\zF^v$ is the vertical lift of a section $\zF$ of $E$, which is
a cocycle $d_\zL\zF=0$ with respect to the exterior derivative of the Lie
algebroid associated with $\zL$ (cf. \cite{IM}). The cocycle property
tells that $[\zF^v,\zL]=0$ and $\zF^v$ is clearly homogeneous, so that
this is precisely the kind of a Jacobi structure described in corollary 1.
This justifies the name Jacobi algebroid given to the bracket on $A(E^*)$
in \cite{GM}. We can slightly generalize the result of \cite{IM} by
considering arbitrary Jacobi structures which are linear with respect to a
vector field $\zD$, i.e. which determine a Lie bracket on linear functions
(homogeneous of degree 1). For functions, vector fields, etc., on a vector
bundle $E$ over $M$, a {\it homogeneous part} is defined. For a function
$f$ let $k$ be the maximal number such that all vertical derivatives of
order $k$ vanish on  $M$  (identified  with  the  0-section).  Then  the
homogeneous part $f_0$ is the
homogeneous polynomial of order $k$ such that all vertical derivatives of
$f-f_0$ of order $(k+1)$ vanish on $M$. For example, if the function
$f$ does not vanish on $M$ then its homogeneous part is just the pull-back
of the function $f$ restricted to $M$. The homogeneous part of a
homogeneous function is just this function. For a vector field $\zG$,
written in local coordinates near the zero-section by
\be
\zG=f_i\pa_{y_i}+g_a\pa_{x^a},
\ee
where $y_i$ are vertical coordinates and $x^a$ are coordinates on the
manifold $M$, the homogeneous part $\zG_0$ of $\zG$ is just the first
non-trivial homogeneous vector field $\zG_0=f'_i\pa_{y_i}$ with $f_i'$
being homogeneous of degree $k$ such that the vertical derivatives of the
vertical coordinates of $\zG-\zG_0$ vanish up to order $(k+1)$.
In particular, if $\zG$ does not vanish on $M$, then
$\zG_0=f_i'\pa_{y_i}$, where $f_i'$ is the  pull-back of $f_i$ restricted
to $M$.

\begin{theorem} Every linear Jacobi structure $J=(\zL',\zG)$ on a
vector bundle $E$ induces a
linear Poisson structure $\zL=\zL'-\zG\we\zD$ such that
$[\zG,\zL]=\zG\we[\zD,\zG]$. It induces also a Jacobi structure
$J'=(\zL+\zG_0\we\zD,\zG_0)$ with $\zG_0$ being the  homogeneous part
of $\zG$.
\end{theorem}
\begin{pf} It is easy to see that the bracket induced by the bivector
$\zL=\zL'-\zG\we\zD$ on linear functions coincides with the Jacobi
bracket, i.e. the bracket is linear and the tensor $\zL$ is homogeneous of
degree -1.  Thus,
\be\label{3}
[\zD,\zL]=[\zD,\zL']-[\zD,\zG]\we\zD=-\zL=-\zL'+\zG\we\zD.
\ee
We get then
\be
[\zL,\zL]=[\zL',\zL']+2[\zG,\zL']-2\zG\we[\zD,\zL']+2\zG\we[\zD,\zG]\we\zD=0
\ee
due to (\ref{3}) and $[\zG,\zL']=0$, $[\zL',\zL']=-2\zG\we\zL'$, so that
$\zL$ is a Poisson tensor. We have additionally
\be\label{4}
[\zG,\zL]=[\zG,\zL']-\zG\we[\zG,\zD]=[\zD,\zG].
\ee
Let now $\zG_0$ be the homogeneous part of $\zG$. For simplicity, assume
that $\zG_0$ is of degree -1 (the general case can be proved in a
completely analogous way). Put $\zG'=\zG-\zG_0$. We have then, according
to (\ref{4}),
\bea
[\zG'+\zG_0,\zL]&=&(\zG'+\zG_0)\we([\zD,\zG']-\zG_0)=\\
&&\zG'\we[\zD,\zG']-\zG'\we\zG_0+\zG_0\we[\zD,\zG'].
\eea
It is easy to see that the right-hand side vanishes on $M$ when applied to
a pair of linear functions. Since the same is true for $[\zG',\zL]$, also
$[\zG_0,\zL]$ vanishes on $M$ when applied to a pair of linear functions.
But $[\zG_0,\zL]$ is a vertical tensor which is constant along fibers, so
that $[\zG_0,\zL]=0$.
\end{pf}

\begin{theorem} Every Jacobi  structure  on  a  vector
bundle $E$ over $M$ which is linear and affine (i.e. such that the linear
and the affine functions are closed with
respect to the Jacobi bracket) is of the form $J=(\zL+\zG\we\zD,\zG)$,
where $\zL$ is a linear Poisson tensor and
$\zG=\zG_0+\zG_1$ is an affine vector field with the decomposition into
homogeneous parts $\zG_0,\zG_1$ of orders -1 and 0, respectively, such
that $[\zG_0,\zL]=0$ and $[\zG_1,\zL]=\zG_0\we\zG_1$. The vertical
vector field $\zG_0$ is the vertical lift $\zF^v$ of certain section
$\zF$ of $E$ which is closed with respect to the exterior derivative
associated with the Lie algebroid structure on $E^*$ induced by $\zL$. If,
additionally, the Jacobi bracket of a linear and a basic function is
basic (i.e. the Jacobi bracket is homogeneous of degree -1), then
$\zG_1=0$ and
\be\label{5}
J=(\zL+\zD\we\zF^v,-\zF^v)
\ee
with $d_\zL\zF=0$.
\end{theorem}
\begin{pf} Since $J=(\zL',\zG)$ is affine, $\zG$  is  an  affine  vector
field splitting into homogeneous parts $\zG=\zG_0+\zG_1$.  According  to
Theorem 1, $\zL=\zL'-\zG\we\zD$ is linear and $[\zG,\zL]=[\zG_0,\zL]+
[\zG_1,\zL]$ equals
\be
(\zG_0+\zG_1)\we[\zD,\zG_0+\zG_1]=\zG_0\we\zG_1.
\ee
Comparing the homogeneous parts of order -1 and 0 we get
$[\zG_0,\zL]=0$ and $[\zG_1,\zL]=\zG_0\we\zG_1$.
\end{pf}
The homogeneous Jacobi structures have been studied in \cite{IM1}.
The bracket corresponding to (\ref{5}) has the form
\be
\{ f,g\}_J=\{ f,g\}+(\zD\we\zF^v)(f,g)-f\zF^v(g)+\zF^v(f)g,
\ee
where $\{\cdot,\cdot\}$ is the linear Poisson  bracket  associated  with
$\zL$.
It is interesting that the bracket
\be
\lna X,Y\rnaf=\lna X,Y\rna+(Deg\we i_\zF)(X,Y) -X\we i_\zF Y+(-1)^xi_\zF
X\we Y,
\ee
has formally the same form with $Deg$ playing the role of the Liouville
vector field and $i_\zf$ being a graded derivative (vector field) of
degree -1 with respect to $Deg$.
We  will  show  latter  that  this  is  not accidental.

\section{Jacobi algebroids and homology}

Recall that for  a  vector  bundle   $E$   over   the   base   manifold
$M$, we denote by  $A(E)=\oplus_{k\in\Z}A^k(E)$ be the exterior algebra
of multisections of
$E$. We will refer to elements of $A^k(E^*)$ as to {\it k-forms} on $E$.
Here we identify $A^0(E)$  with  the  algebra  $C^\infty(M)$  of  smooth
functions on the base and $A^k(E)=\{ 0\}$ for $k<0$.
Denote by $\vert X\vert$ the Grassmann degree of the multisection $X\in
A(E)$. We  will use  the convention that we write $x$ for $\vert
X\vert-1$ -- the  shifted degree  of $X$ (this is the Lie algebra degree
of $X$ with respect to  the Schouten bracket $\lna\cdot,\cdot\rna$
induced  by  any  Lie  algebroid bracket on $E$).

\medskip
\medskip
This is the idea going back to E.~Witten \cite{Wi} to deform the de Rham
exterior derivative by gauging the cotangent bundle by the multiplication
by the function $e^f$:
\be
\D^{\D f}\zm=e^{-f}\D(e^f\zm)=\D\zm+\D f\we\zm.
\ee
We have clearly $(\D^{\D f})^2=0$ and we get the corresponding
cohomology being equivalent to de Rham cohomology. This time, however,
$\D^{\D f}$ is not a derivation but a first-order differential
operator with respect to the wedge product on differential forms.

\medskip
A natural generalization is to start with the exterior derivative $\D$
associated with a Lie algebroid structure in a vector bundle $E$ over
$M$ and to take any 1-cocycle $\zF$ instead of the coboundary $\D f$, so
that $\D^\zF\zm=\D\zm+\zF\we\zm$. This is exactly the exterior
differential we obtain
for a Jacobi algebroid associated with the 1-cocycle $\zF$ by an analog of
the Cartan formula (cf. \cite{IM, GM}):
\bea\nn
        \D^\zF\zm(X_1,\dots,X_{k+1}) &=& \sum_i (-1)^{i+1}
\lna X_i,\zm(X_1,\dots,\widehat{X}_i,\dots ,X_{k+1})\rnaf \\
&+& \sum_{i<j} (-1)^{i+j}\zm (\lna X_i,X_j\rnaf, X_1,\dots , \widehat{X}_i,
\dots ,\widehat{X}_j, \dots , X_{k+1}).\label{D}
\eea
Even if the 1-cocycle $\zF$ is not
exact, there is a nice construction \cite{IM} which allows to view $\zF$
as being exact but for an extended Lie algebroid in the bundle $\hat
E=E\ti\R$ over $M\ti\R$. The sections of this bundle may be viewed as
time-dependent sections of $E$. The sections of $E$ form a Lie subalgebra
of time-independent  sections  in  the  Lie  algebroid  $\hat  E$  which
generate the $C^\infty(M\ti\R)$-module of sections of $\hat E$  and  the
whole structure is uniquely determined by putting the anchor $\hat\zr(X)$
of a time-independent section $X$ to be $\hat\zr(X)=\zr(X)+\langle\zF,X
\rangle\pa_t$,
where $t$ is the standard coordinate function in $\R$ and $\zr$ is
the anchor in $E$. All this is consistent (thanks to the fact that $\zF$
is a cocycle) and defines a Lie algebroid structure on $\hat E$ with the
exterior derivative $\D$ satisfying $\D t=\zF$.

\medskip
Let now $U:A(E)\ra A(\hat E)$ be natural embedding of the Grassmann algebra
of $E$ into the Grassmann subalgebra of time-independent sections of $\hat
E$. It is obvious that $U$ is a homomorphism of the corresponding Schouten
brackets:
\be
\lna U(X),U(Y)\rnah=U(\lna X,Y\rna),
\ee
where  we  use  the  notation  $\lna\cdot,\cdot\rna$   and   $\lna\cdot,
\cdot\rnah$ for the Schouten brackets in $E$ and $\hat E$, respectively.
Let us now gauge $A(E)$ inside $A(\hat E)$ by putting
\be
\tilde U(X)=e^{-xt}U(X)
\ee
for any homogeneous element $X$. Note that $\tU$ preserves  the  grading
but not the wedge product. It can be easily proved (cf. \cite{GM} where
this is proved for an extension of $\hat E$) that the Jacobi  algebroid
bracket  (\ref{jb})  can  be obtained by this gauging $A(E)$ in $A(\hat E)$.
\begin{theorem} For any homogeneous elements $X,Y\in A(E)$ we
have
\be\label{j}
\lna\tU (X),\tU(Y)\rnah=\tU(\lna X,Y\rna+xX\we
i_\zF Y-(-1)^xyi_\zF X\we Y)=\tU(\lna X,Y\rnaf).
\ee
\end{theorem}
From the above theorem we get for free the following.
\begin{corollary}(\cite{IM}) The Schouten-Jacobi bracket
(\ref{jb}) is a graded Lie bracket for $A(E)$.
\end{corollary}
Thus we have obtained the Jacobi algebroid bracket and the corresponding
exterior differential by gauging. What about the other ingredients of the
Cartan calculus? The contraction is obvious, so let us define the Lie
differential. For a Lie algebroid $E$ the Lie differential $\Ll_X$ along
a multisection $X$ of $E$ acting on $A(E^*)$ is defined by
\be
\Ll_X\zm=i_X\D\zm+(-1)^x\D i_X\zm.
\ee
We have the following well-known formulae (cf. e.g. \cite{KSM,Mi})
\bea\label{ld}
\lbo\Ll_X,\Ll_Y\rbo&=&-\Ll_{\lna Y,X\rna},\\
\lbo\Ll_X,i_Y\rbo&=&-i_{\lna Y,X\rna},
\eea
where
\bea
\lbo\Ll_X,\Ll_Y\rbo&=&\Li_X\circ\Li_Y-(-1)^{xy}\Li_Y\circ\Li_X,\\
\lbo\Ll_X,i_Y\rbo&=&\Li_X\circ i_Y-(-1)^{(y+1)x}i_Y\circ\Li_X,
\eea
are the graded commutators of the graded morphism of $A(E^*)$.
If we define the deformed Lie differential by
\be
\Llf_X\zm=i_X\D^\zF\zm+(-1)^x\D^\zF i_X\zm,
\ee
then obviously
\be
\lbo\Llf_X,\Df\rbo=\Llf_X\circ\Df-(-1)^x\Df\circ\Llf_X=0.
\ee
Using the formula
\be
i_X(\zF\we\zm)=i_{X_\zF}\zm-(-1)^x\zF\we i_X\zm,
\ee
where we write $X_\zF$ for $i_\zF X$, we get
\be\label{ld1}
\Llf_X=\Ll_X+i_{X_\zF}.
\ee
Note that this coincides with the definition of $\zF$-Lie derivative  in
\cite{IM} for $X$ being just sections of $E$. In this case
$i_{X_\zF}\zm=\langle X,\zF\rangle\zm$.
However, in spite of the fact that the Lie differential was deformed,  we
get the same formulae as (\ref{ld}) with the original  Schouten  bracket
instead of the Schouten-Jacobi bracket.
\begin{theorem} The following identities hold:
\bea\label{a}
\lbo\Llf_X,\Llf_Y\rbo&=&-\Llf_{\lna Y,X\rna},\\
\lbo\Llf_X,i_Y\rbo&=&-i_{\lna Y,X\rna}.\label{b}
\eea
\end{theorem}
\begin{pf} Since $\lbo i_X,i_Y\rbo=0$, we get
\be
\lbo\Llf_X,i_Y\rbo=\lbo\Ll_X+i_{X_\zF},i_Y\rbo=-i_{\lna Y,X\rna}.
\ee
Now,
\be
\lbo\Llf_X,\Llf_Y\rbo=\lbo\Llf_X,\lbo i_Y,\Df\rbo\rbo=
-\lbo i_{\lna Y,X\rna},\Df\rbo=-\Llf_{\lna Y,X\rna},
\ee
since $\Llf_X=\lbo i_X,\Df\rbo$ and $\lbo\Llf_X,\Df\rbo=0$.
\end{pf}

\medskip\noindent
To obtain similar formulae but with the Schouten-Jacobi bracket we will
use again a proper gauging. Let us use the embedding $\hat U$ of the
Grassmann algebra $A(E^*)$ into $A(\widehat{E^*})=A((\widehat{E})^*)$ by
$\hat U_a(\zm)=e^{(\vm+a)t}U(\zm)$.
Here the elements $U(\zm)$ are regarded as time-independent
sections. Using the standard Lie differential $\Li$ and the exterior
derivative $\D$ for the extended  Lie  algebroid  $\hat  E$, we get the
following.
\begin{lemma}
\bea\label{li1}
i_{\tU(X)}\hU_a(\zm)&=&e^t\hU_a(i_X\zm);\\
\Li_{\tilde U(X)}\hat U_a(\zm)&=&\hat U_a(\Li_X\zm+(\vm+a) i_{X_\zF}\zm
-(-1)^xx\zF\we i_X\zm).\label{li2}
\eea
\end{lemma}
\begin{pf} We will write, for simplicity, $X$ and $\zm$ instead of  $U(X)$
and  $U(\zm)$,  regarding  $A(E)$  and  $A(E^*)$   as   time-independent
sections in $A(\widehat E)$ and $A(\widehat{E^*})$, respectively. The
Schouten-Nijenhuis bracket  on
embedded $A(E)$ and the exterior derivative on $A(E^*)$ coincide with the
restrictions of the Schouten bracket on $A(\hat  E)$  and  the  exterior
derivative on $A(\hat E^*)$, respectively. We have clearly
\be
i_{\tU(X)}\hU_a(\zm)=e^{-xt}i_X(e^{(\vm+a)t}\zm)=e^t\hU_a(i_X\zm).
\ee
Using the fact that $\D t=\zF$, we get in turn
\bea
\Li_{\tU(X)}\hat U_a(\zm)&=&i_{\tilde U(X)}\D(e^{(\vm+a) t}\zm)+(-1)^x\D
i_{\tilde U(X)}(e^{(\vm+a) t}\zm)=\\
&&e^{-xt}i_X((\vm+a)e^{(\vm+a) t}\zF\we\zm+e^{(\vm+a) t}\D\zm)
+(-1)^x\D(e^{(\vm+a-x)t}i_X\zm)=\nn \\
&&\hU_a(\Li_X\zm+(\vm+a) i_{X_\zF}\zm-(-1)^xx\zF\we i_X\zm).\nn
\eea
\end{pf}

\medskip\noindent
For the Jacobi algebroid associated with the 1-cocycle $\zF$ we will call
\be
\Llh^{\zF,a}_X\zm=\Li_X\zm+(\vm+a)i_{X_\zF}\zm-(-1)^xx\zF\we i_x\zm
\ee
the {\it Jacobi-Lie differential} of $\zm$ along $X$. According to the
above lemma,
\be\label{ld2}
\Li_{\tilde U(X)}\hat U_a(\zm)=\hat U_a(\Llh^{\zF,a}_X\zm)
\ee
and we have the following analog of (\ref{a}).

\begin{theorem}
\be\label{a'}
\lbo\Llh^{\zF,a}_X,\Llh^{\zF,a}_Y\rbo=-\Llh^{\zF,a}_{\lna Y,X\rnaf}.
\ee
\end{theorem}
\begin{pf} In view of (\ref{ld2}),
\be
\hU_a(\Llh^{\zF,a}_X\circ\Llh^{\zF,a}_Y(\zm))=\Li_{\tU(X)}\circ
\Li_{\tU(Y)}\hU_a(\zm),
\ee
so that
\bea
\hU_a(\lbo\Llh^{\zF,a}_X,\Llh^{\zF,a}_Y\rbo(\zm))&=&
\lbo\Li_{\tU(X)},\Li_{\tU(Y)}\rbo\hU_a(\zm)=
-\Li_{\lna\tU(Y),\tU(X)\rnah}\,\hU_a(\zm)=\\
&&-\Li_{\tU(\lna Y, X\rnaf)}\hU_a(\zm)=
\hU_a(-\Llh^{\zF,a}_{\lna Y,X\rnaf}(\zm))
\eea
and theorem follows, since $\hU$ is injective.
\end{pf}

\medskip\noindent
Instead of (\ref{b}) we have the following deformed version.
\begin{theorem}
\be\label{ab'}
\lbo\Llh^{\zF,a}_X,i_Y\rbo=-i_{\lna Y,X\rnaf}-(-1)^xi_{Y\we X_\zF}.
\ee
\end{theorem}
\begin{pf}
Analogously as above we get by (\ref{li1})
\be\label{ab}
\lbo\Li_{\tU(X)},e^{-t}i_{\tU(Y)}\rbo \hU_a(\zm)=
\hU_a(\lbo\Llh^{\zF,a}_X,i_Y\rbo(\zm)).
\ee
Now, using
\be
\Li_Z(e^{-t}\zm)=e^{-t}(\Li_Z\zm+(-1)^zi_{Z_\zF}\zm),
\ee
we get
\bea
\hU_a(\lbo\Llh^{\zF,a}_X,i_Y\rbo(\zm))&=&
e^{-t}(\lbo\Li_{\tU(X)},i_{\tU(Y)}\rbo+(-1)^xi_{\tU(X)_\zF}\circ
i_{\tU(Y)})\hU_a(\zm)=\\
&&e^{-t}(-i_{\lna\tU(Y),\tU(X)\rna\hat{}}+(-1)^xe^{-t}i_{\tU(X_\zF)} \circ
i_{\tU(Y)})\hU_a(\zm)=\\
&&\hU_a(-i_{\lna Y,X\rnaf}-(-1)^xi_{Y\we X_\zF})
\eea
and the theorem follows.
\end{pf}

\medskip\noindent
An element $P\in A^2(E)$ with $\lna P,P\rnaf=0$  we  will  call  a  {\it
Jacobi element} (or {\it canonical structure}).
An immediate consequence of (\ref{li2}) is the following observation.

\begin{corollary} For a Jacobi element $P$ of $A(E)$ the Lie
differential
\be\label{LJ}
\Llh^{\zF,a}_P\zm=\Li_P\zm+(\vm+a)i_{P_\zF}+\zF\we i_P\zm
\ee
is a homology operator on  $A(E^*)$, i.e. $\vert\Llh^{\zF,a}_P\zm\vert
=\vm-1$  and  $(\Llh^{\zF,a}_P)^2=0$.  Moreover,  $\Llh^{\zF,a}_P$  is   a
generating operator for the Schouten-Nijenhuis bracket on $A(E^*)$:
\be\label{Kz}
\lna\zm,\zn\rna_P=(-1)^\vm(\Llh^{\zF,a}_P(\zm\we\zn)-\Llh^{\zF,a}_P
(\zm)\we\zn-(-1)^\vm\zm\we\Llh^{\zF,a}_P(\zn))
\ee
which does not depend on $a$ and which is the Schouten-Nijenhuis bracket
of the Lie algebroid bracket on $E^*$, defined for $\zm,\zn\in A^1(E^*)$
by
\bea\label{Koz1}
\lna\zm,\zn\rna_P&=&i_{P_\zm}\Df\zn-i_{P_\zn}\Df\zm+\Df\langle P,
\zm\we\zn\rangle=\\
&&\Llf_{P_\zm}\zn-\Llf_{P_\zn}\zm-\Df\langle P,\zm\we\zn\rangle.
\eea
\end{corollary}

\medskip\noindent
Since $\Llh^{\zF,a}_P-\Llh^{\zF,a'}_P=(a-a')i_{P_\zF}$ is a  derivation,
all   Lie   differentials   $\Llh^{\zF,a}_P$   are   equally   good   as
generators of the Schouten-Nijenhuis bracket (\ref{LJ}). In \cite{ILMP},
Theorem  4.8,  the  authors   have   introduced   $\Llh^{\zF,0}_P$   and
$\Llh^{\zF,1}_P$. The first Lie differential $\Llh^{\zF,0}_X$ (denoted
simply  $\Llh^{\zF}_X$) reduces to $\Llf_X$ on 1-forms for $X$
being  just  a   section   of   $E$.   The   second   Lie   differential
$\Llh^{\zF,1}_X$, applied for functions, describes, for $X\in A^1(E)$,
the bracket: $\Llh^{\zF,1}_X(f)=\lna X,f\rnaf$. The homology defined by
$\Llh^{\zF}_P$   for   a   Jacobi   element    $P$    we    will    call
{Lichnerowicz-Jacobi  homology}  and  denote  by  $H_*^{LJ}(E,P)$.   The
homology defined by $\Llh^{\zF,1}_P$ we will  call  {Jacobi homology}
and denote by $H_*^{J}(E,P)$.

The fact that $(\Llh^{\zF,a}_P)^2=0$ for all $a$ has the obvious
consequence that $\Llh^{\zF}_P$ and $i_{P_\zF}$ commute.

\begin{theorem} For a Jacobi element $P$
\be\label{comm}
\lbo\Llh^{\zF}_P,i_{P_\zF}\rbo=\Llh^{\zF}_P\circ i_{P_\zF}+
i_{P_\zF}\circ\Llh^{\zF}_P=0.
\ee
\end{theorem}

\medskip\noindent
{\bf Remark.} The bracket (\ref{Koz1}) has been introduced  in  \cite{IM}
as a generalization of the  triangular  Lie  bialgebroid  \cite{MX}  for
Jacobi algebroids (generalized Lie algebroids)  and  it  is  an  obvious
generalization of the Koszul-Fuchssteiner  \cite{Kz,Fu}  bracket  on
1-forms induced by a Poisson structure.

\medskip\noindent
{\bf Example 1.} If $E=TM\oplus\R$ is  a
Lie algebroid of first-order differential  operators  on  $C^\infty(M)$,
i.e. the Lie bracket on sections reads
\be
\lna(X,f),(Y,g)\rna=([X,Y],X(g)-Y(f)),
\ee
and the 1-cocycle is $\zF((X,f))=f$, then every Jacobi  element  $P=(\zL,
\zG)$ being a section of $A^2(E)={\bigwedge}^2TM\oplus  TM$  is  just  a
Jacobi structure on $M$. Identifying elements of $\zW^k(E)$  with  pairs
$(\zm,\zn)$,  where  $\zm$  is  a  $k$-form  on  $M$  and  $\zn$  is   a
$(k-1)$-form on $M$, we get $\D(\zm,\zn)=(\D\zm,-\D\zn)$ and
\be
\Li_{(\zL,\zG)}(\zm,\zn)=(\Li_\zL\zm-\Li_\zG\zn,-\Li_\zL\zn).
\ee
Thus,
\be
\Llh_P^\zF(\zm,\zn)=(\Li_\zL\zm-\Li_\zG\zn+ki_\zG\zm,-\Li_\zL\zn+
i_\zL\zm-(k-1)i_\zG\zn),
\ee
so we get exactly the standard Lichnerowicz-Jacobi homology operator
as described in  \cite{LLMP,ILMP,Va1,BMMP}.
For a Poisson structure it  reduces  to  the  Koszul-Brylinski  boundary
operator (see e.g. \cite{Br} or \cite{Va}, Ch.5).
The    homology     operator     $\Llh^{\zF,1}_P$ is related  to
Chevalley-Eilenberg homology operator $\zd^{CE}$  of  the  corresponding
Jacobi bracket with respect to the maps
\be
\zp_k:C^\infty(M)\ti\zL^kC^\infty(M)\ra\zW^k(E),\quad
\zp_k(f\ti f_1\we\dots\we f_k)=f\D^\zF f_1\we\dots\we\D^\zF f_k
\ee
by $\Llh_P^{\zF,1}\circ\zp_k=\zp_{k-1}\circ\zd^{CE}$.

\section{Courant-Jacobi algebroids}

The  method  of  gauging  can  be  used  to  deform   the   {\it Courant
brackets} associated with a Lie algebroid bracket  $[\cdot,\cdot]$  on
$E$. Recall that the Courant bracket \cite{Co} is the following  bracket
on $E\oplus E^*$:
\be\label{cb}
[X+\zx,Y+\zh]_{C}=[X,Y]+\Li_X\zh-\Li_Y\zx+\frac{1}{2}\D(i_Y\zx-i_X\zh).
\ee
{\it  Dirac  structures},  introduced  independently  by   Courant   and
Weinstein \cite{CW} and Dorfman \cite{Do1} for $E=TM$, can be defined  as
subbundles $L$ of $E\oplus E^*$ which are maximally isotropic under the
canonical symmetric pairing $\la X+\zx,Y+\zh\ran_+=i_Y\zx+i_X\zh$ and
closed  with  respect   to   the  Courant  bracket.  For  their  use  in
studying completely integrable systems of partial differential equations
we refer to \cite{Do2}.

Now, choosing a 1-cocycle $\zF\in\zW^1(E)$, we can consider the extended
Lie algebroid $\hat E$ defined in the previous section and the embedding
\be
\U(X+\zx)=U(X)+\hat U_0(\zx)
\ee
of   $\textrm{Sec}(E\oplus   E^*)$   into   $\textrm{Sec}(\hat   E\oplus
\hat{E}^*)$. In other words, $\U(X+\zx)=X+e^t\zx$  when  $X,\zx$  are
regarded as time-independent sections of $\hat E$  and  $\hat{E}^*$,
respectively.
For computational aims we can just think that this is a
true gauging, i.e. we work in $E$ and $E^*$ and  $t\in  C^\infty(M)$  is
the potential for $\zF$, i.e. $\D t=\zF$.

In $\hat E\oplus \hat E^*=\hat  E\oplus\widehat{E^*}$  we have its own
Courant bracket $[\cdot,\cdot]_{C}^{\wedge}$. Now, following the ideas of the
previous  sections,  one  proves  easily that
\be
[\U(X+\zx),\U(Y+\zh)]^{\wedge}_{C}=\U([X,Y]+\Li^\zF_X\zh-\Li^\zF_Y\zx+\frac{1}
{2}\D^\zF(i_Y\zx-i_X\zh)).
\ee
The bracket
\be\label{cjb}
[X+\zx,Y+\zh]^\zF_{C}=
[X,Y]+\Li^\zF_X\zh-\Li^\zF_Y\zx+\frac{1}{2}\D^\zF(i_Y\zx-i_X\zh)
\ee
we will call the {\it Courant-Jacobi bracket} (associated  with  $\zF$).
Similarly as above, a {\it Dirac-Jacobi structure} is a subbundle $L$ of
$E\oplus  E^*$ which is maximally isotropic under the canonical
symmetric  pairing   and  closed  with  respect  to  the  Courant-Jacobi
bracket.

\medskip\noindent
{\bf Example 2.} Let $E=TM\oplus\R$ be the Lie algebroid  of  first-order
differential operators on $M$. The sections of $E$ are  identified  with
pairs  $(X,f)$,  where  $X$  is a vector  field   on   $M$   and  $f\in
C^\infty(M)$, and the Lie algebroid bracket reads
\be
[(X_1,f_1),(X_2,f_2)]=([X_1,X_2],X_1(f_2)-X_2(f_1)).
\ee
The sections of the dual bundle are pairs $(\za,h)$, where  $\za$  is  a
1-form on $M$ and $h\in C^\infty(M)$, with the obvious pairing. The  Lie
algebroid  exterior  derivative  of  a  function  is  $\D  f=(\D  f,0)$,
where  $\D$  on  the  right-hand  side   is   the   standard   exterior
derivative  (we  hope  that  this  abuse  in  notation  will  cause   no
confusion) and the Lie derivative is given by
\be
\Li_{(X,f)}(\za,h)=(\Li_X\za+h\D f,X(h)).
\ee
Thus the Courant bracket on $E\oplus E^*$ reads
\bea\nn
&&[((X_1,f_1)+(\za_1,g_1),(X_2,f_2)+(\za_2,g_2)]_{C}=
([X_1,X_2],X_1(f_2)-X_2(f_1))\\
&&+(\Li_{X_1}\za_2-\Li_{X_2}\za_1+g_2\D f_1-g_1\D f_2
+\frac{1}{2}\D(
i_{X_2}\za_1-i_{X_1}\za_2+f_2g_1-f_1g_2),X_1(g_2)-X_2(g_1)).
\eea
The standard 1-cocycle $\zF$  on  $E$  leading  to  the  Schouten-Jacobi
algebra  of first-order   polydifferential   operators  is   given   by
$\zF((X,f))=f$.  It  is  easy  to  see  now   that   the   corresponding
Courant-Jacobi bracket is given by
\bea\nn
&&[((X_1,f_1)+(\za_1,g_1),(X_2,f_2)+(\za_2,g_2)]^\zF_{C}=
([X_1,X_2],X_1(f_2)-X_2(f_1))\\
&&+(\Li_{X_1}\za_2-\Li_{X_2}\za_1+g_2\D f_1-g_1\D f_2
+\frac{1}{2}\D(
i_{X_2}\za_1-i_{X_1}\za_2+f_2g_1-f_1g_2)+f_1\za_2- f_2\za_1,
\\ \nn
&&X_1(g_2)-X_2(g_1)+\frac{1}{2}(i_{X_2}\za_1-
i_{X_1}\za_2-f_2g_1+f_1g_2)).
\eea
This is exactly the bracket introduced by A.~Wade  \cite{Wa}  to  define
${\cal E}^1(M)$-Dirac structures.

\medskip
The  Courant  bracket  (\ref{cb})  can  be  generalized  to  a   bracket
associated with a Lie bialgebroid
\be
((E,[\cdot,\cdot]_E),(E^*,[\cdot,\cdot]_{E^*}))
\ee
in  the sense  of Mackenzie and  Xu \cite{MX} (which, in turn, is a
fundamental example of  a  {\it  Courant algebroid} bracket
\cite{LWX}):
\bea\label{cab}
[X+\zx,Y+\zh]_{C}&=&([X,Y]_E+\Li^{E^*}_{\zx}Y-\Li^{E^*}_{\zh}X
-\frac{1}{2}\D_{E^*}(i_Y\zx-i_X\zh)\\
&&+([\zx,\zh]_{E^*}+\Li^E_X\zh-\Li^E_Y\zx+\frac{1}{2}\D_E(i_Y\zx-
i_X\zh)).
\eea

Consider now a Jacobi bialgebroid \cite{GM} (generalized Lie bialgebroid
in the sense of \cite{IM})
\be\label{jbal}
((E,[\cdot,\cdot]^{\zF_0}_E),(E^*,[\cdot,\cdot]^{X_0}_{E^*}))
\ee
associated with 1-cocycles $\zF_0\in\textrm{Sec}(E^*)$,  $X_0\in\textrm{
Sec}(E)$ with respect to the Lie algebroid brackets on  $E$  and  $E^*$,
respectively.                                                        Let
$(\hat E,[\cdot,\cdot]^{\wedge}_E)$ and $(\hat{E}^*, [\cdot,\cdot]
^{\wedge} _{E^*})$ be extended Lie algebroids associated with these
cocycles  with
$\D_E(t)=\zF_0$ and $\D_{E^*}(t)=X_0$.
Let $[\cdot,\cdot] ^{{0}} _{E^*}$ be the  Lie  algebroid  bracket  on
$\hat{E}^*$ obtained from $[\cdot,\cdot] ^{\wedge} _{E^*}$ by gauging
by $e^{-t}$:
\be
[\zx,\zh] ^{{0}}_{E^*}=e^t[e^{-t}\zx,e^{-t}\zh] ^{\wedge} _{E^*}.
\ee
It has been proved in \cite{IM}, Theorem 4.11, that
\be
((\hat E,[\cdot,\cdot]^{\wedge}_E),(\hat{E}^*, [\cdot,\cdot]
^{{0}}_{E^*}))
\ee
is a Lie bialgebroid, so  we  can  consider  the  corresponding  Courant
algebroid bracket $[\cdot,\cdot]_{C}^{\wedge}$ (\ref{cab}) on $\hat
E\oplus\hat{E}^*$. Now, using
$\U:E\oplus E^*\ra \hat E\oplus\hat{E}^*$, we get
\be
[\U(X+\zx),\U(Y+\zh)]^{\wedge}_{C}=\U([X+\zx,Y+\zh]_{C}^{\zF_0,X_0}),
\ee
where
\bea\label{cjab}
[X+\zx,Y+\zh]^{\zF_0,X_0}_{C}&=&([X,Y]_E+\Li^{E^*}_{\zx}Y-
\Li^{E^*}_{\zh}X\\ \nn
&&-\frac{1}{2}\D_{E^*}(i_Y\zx-i_X\zh)+\frac{1}{2}(i_Y\zx-i_X\zh)X_0)\\ \nn
&&+(([\zx,\zh]_{E^*}+{(\Li^E)}^{\zF_0}_X\zh-{(\Li^E)}^{\zF_0}_Y\zx+
\frac{1}{2}\D^{\zF_0}_E(i_Y\zx-i_X\zh)).
\eea
This is the {\it Courant-Jacobi  bracket}  associated  with  the  Jacobi
bialgebroid. In the  case  when  $[\cdot,\cdot]_{E^*}$  is  trivial  and
$X_0=0$ we end up with (\ref{cjb}).

Being obtained  by  gauging this bracket has  properties similar to
that  of Courant algebroid bracket. The  abstract  of  these  properties
leads to the following definition (based on  the  definition  of  Courant
algebroid proposed in \cite{Ro} rather than on the  original  definition
in \cite{LWX}). We have reduced the number of axioms using ideas similar
to \cite{Uch}.

\begin{definition} A {\sl Courant-Jacobi algebroid} is a  vector  bundle
$F$ over $M$ together  with
\begin{description}
\item{(i)} a  nondegenerate  symmetric  bilinear  form
$\la\cdot,\cdot\ran$ on the bundle,
\item{(ii)} a Loday  operation  $\circ$  on $\textrm{Sec}(F)$  (i.e.  a
bilinear  operation  satisfying  the  Jacobi identity; the terminology
goes back  to  \cite{KS1}),
\item{(iii)} a bundle map $\zr:F\ra TM\oplus\R$ which is a  homomorphism
into the Lie algebroid of first-order differential operators
\be
[\zr(e_1),\zr(e_2)]=\zr(e_1\circ e_2),
\ee
\end{description}
satisfying the following properties:
\begin{enumerate}
\item $\la e_1\circ e,e\ran=\la e_1,e\circ e\ran$,
\item $\zr(e_1)(\la e,e\ran)=2\la e_1\circ e,e\ran$,
\end{enumerate}
for all $e,e_1,e_2\in\textrm{Sec}(F)$.
\end{definition}

Note that we can reformulate  the  above  definition  in  terms  of  the
first-order differential operator
\be
{\cal D}:C^\infty(M)\ra Sec(F), \quad \la{\cal D}(f),e\ran=\zr(e)(f).
\ee
The  Courant-Jacobi  algebroid  (like  the  Courant  algebroid)   is   a
particular case of a pure algebraic structure  described  in  \cite{JL},
Theorem 3.3.
For a Jacobi bialgebroid (\ref{jbal}) the corresponding Courant-Jacobi
algebroid bracket on $F=E\oplus E^*$ is
\bea\label{cjalb}
(X+\zx)\circ(Y+\zh)&=&([X,Y]_E+\Li^{E^*}_{\zx}Y-i_\zh\D_{E^*}X
+\la Y,\zx\ran X_0)\\ \nn
&&+([\zx,\zh]_{E^*}+{(\Li^E)}^{\zF_0}_X\zh-i_Y\D^{\zF_0}_E\zx).
\eea
The symmetric form is $\la\cdot,\cdot\ran_+$, and $\zr=\zr_E+\zr_{E^*}
+i_{\zF_0+X_0}$.

The  definition  of  a  {\it  Dirac   structure}   associated   with   a
Courant-Jacobi algebroid is obvious.
We will postpone a study of  Courant-Jacobi  algebroids  to  a  separate
paper. Note only the following.

\begin{theorem}  Every   Dirac   structure   $L$   associated   with   a
Courant-Jacobi algebroid structure  on  $E$  induces  on  $L$  a  Jacobi
algebroid structure with the bracket being the restriction of $\circ$ to
sections  of  $L$  and  the   1-cocycle   being   the   restriction   of
$\la\zF,\cdot\ran$ to sections of $L$.
\end{theorem}

\section{Krasil'shchik calculus for first-order  differential operators}

It this section we will present the Krasil'shchik's approach to Schouten
brackets \cite{Kr1,Kr2} adapted to Jacobi structures, i.e. we will  show
that the identification of the concepts of Schouten and Poisson brackets
can be extended  to  canonical  brackets  of  first-order
polydifferential operators  and  Jacobi   brackets.

The   differential   calculus   for associative commutative algebras has
been   developed    by    A.~M.~Vinogradov    \cite{ViA}    (see    also
\cite{KLV,VK,ViM}).

In \cite{GM} we have observed that the Schouten bracket  on  multivector
fields (i.e. elements of $A(TM)$) can be viewed as  the  restriction  of
the   Richardson-Nijenhuis   bracket   on   multilinear   operators   on
$C^\infty(M)$ to polyderivations. At the same time, when  reducing  the
Richardson-Nijenhuis bracket to first-order polydifferential operators,
we get  not  a  Schouten-type  but  a  Jacobi-type  bracket  which is a
particular case of what was recently studied  in  \cite{IM}  under  the
name of a generalized Lie algebroid.
Since  a   very   general   approach   to    Schouten    brackets    and
supercanonical
structures has been already developed  in  \cite{Kr2},  we  will  follow
these ideas to develop a similar calculus  `on  the  Jacobi  level'.  We
slightly  change  the  approach  of  \cite{Kr2}  using  rather  a  shift
$\za\in\Z^n$ in the original grading than a divided grading. Of  course,
any such shift defines a divided grading by parity  of  coefficients  of
$\za$ so that the signs in the Krasilshchik's and our approaches  remain
the same,  but using  the  shift  allows  to  trace  better  the  graded
structures which are introduced and fits better to the concept of graded
Poisson or Jacobi algebras. We will mostly skip  the  proofs  which  are
standard inductions and just matters of simple  calculations  completely
parallel to those in \cite{Kr2}. All elements considered in formulae are
uniform, i.e. with a well-defined degree.

We start with an $n$-graded associative  commutative  algebra  $\A$  with
unity  $\1$  and  we  will  write  simply   $a$   instead   of    $\vert
a\vert\in\Z^n$.
Let us fix $\za\in\Z^n$.

We  define   the  graded  $\A$-bimodules  $\Di_i(\A)$  (denoted  shortly
by $\Di_i$ if $\A$ is fixed)  of  genus $i$ of
polydifferential operators of first-order by induction.  Note  that  the
genus relates the right-module structure  to  the  left-module
structure by
\be\label{genus}
p\cdot a=(-1)^{\la a,p+i\za\ran}a\cdot p
\ee
for $p\in\Di_i$, $a\in\A$, so that we get the right-module structure from
the left one by definition.
For $\Di_0$ we take just $\A$ (we can start with an arbitrary module  of
genus 0, but this choice is sufficient for our purposes in this paper).
Then, we  take  $\Di_1$  as  the  space  of  those  linear  graded  maps
$D:\A\ra\A$ which satisfy
\be
D(ab)=D(a)b+(-1)^{\la a,D\ran}aD(b)-D(\1)ab,
\ee
i.e. $\Di_1$ is the module  of  first-order  differential  operators  on
$\A$. For $i>1$  we  define  inductively  $\Di_i$  as  formed  by  those
graded linear maps $D:\A\ra\Di_{i-1}$ for which
\bea
D(ab)&=&D(a)\cdot b+(-1)^{\la a,D+(i-1)\za\ran}
a\cdot D(b)-D(\1)\cdot ab,\\
D(a,b)&=&-(-1)^{\la a+\za,b+\za\ran}D(b,a).
\eea
The notation is clearly
\be
D(a_1,\dots,a_j)=D(a_1,\dots,a_{j-1})(a_j).
\ee
The left module structure on $\Di_i$ is obvious:  $(a\cdot  D)(b)=a\cdot
D(b)$, the right one is determined by genus.
There is a natural graded subspace
\be
\De(\A)=\oplus_{i=0}^\infty\De_i(\A)
\ee
of polyderivations, i.e. $\De_i$ consists of  those  elements  $D$  from
$\Di_i$  such  that  $D(\1)=0$  (by  graded  symmetry  this  means  that
$D(\cdots,\1,\cdots)=0$).

As in the classical case,
$\Di_1=\De_1\oplus\A\cdot I$, where $I\in\Di_1$ is the identity operator
on $\A$.  We  will  show  a  generalization  of  this  fact  later  (see
(\ref{dec})).
We can extend the product in the algebra $\A$ to the space
\be
\Di=\oplus_{i=0}^\infty\Di_i
\ee
of first-order polydifferential operators putting inductively
\be\label{prod}
(A\cdot B)(a)=(-1)^{\la a+\za,B+j\za\ran+j}A(a)\cdot B+A\cdot B(a)
\ee
for $a\in\A,A\in\Di_i,  B\in\Di_j$.  Similarly  as  in  \cite{Kr2},  one
checks that $A\cdot B$ is in $\Di_{i+j}$ and one proves the following.

\begin{theorem}  The  multiplication  (\ref{prod})   turns   the   space
$\Di(\A)$ of first-order polydifferential operators on $\A$ into an
$(n+1)$-graded associative commutative unital  algebra  with  homogeneous
elements $A$  from $\Di_i$ being  of  degree   $(\vert  A\vert+i\za,i)$.
The graded subspace $\De(\A)$ is a graded subalgebra of $\Di(\A)$.
\end{theorem}

\medskip\noindent
{\bf Remark.} The graded commutativity  in  the  above  theorem  can  be
explicitly written in the form
\be
A\cdot B=(-1)^{\la A+i\za,B+j\za\ran +ij}B\cdot A.
\ee
Here and later on we denote by $A$ also the degree $\vert A\vert$ if  no
confusion arises.
Note also that for $A,B\in\Di_1$ we get
\be
(A\cdot B)(a,b)=(-1)^{\la  a+B,A+\za\ran}B(a)A(b)-(-1)^{\la  a+\za,B+\za
\ran}A(a)B(b).
\ee
In particular, for the ungraded case,
\be
A\cdot B(a,b)=B(a)A(b)-A(a)B(b)=(B\we A)(a,b).
\ee
Thus the defined product is the reversed wedge product.

\medskip
The graded Schouten-Jacobi bracket on $\Di$ we  define  formally  as  in
\cite{Kr2} putting $\lna a,b\rna=0$ for $a,b\in\A$,
\be\label{bra}
\lna D,a\rna=D(a),\quad \lna a,D\rna=(-1)^{\la a-\za,D+(i-1)\za\ran+i}
D(a),
\ee
for $a\in\A$, $D\in\Di_i$, and
\be
\lna A,B\rna(a)=(-1)^{\la a-\za,B+(j-1)\za\ran+j-1}\lna A(a),B\rna+
\lna A,B(a)\rna
\ee
for  $a\in\A,A\in\Di_i,  B\in\Di_j$.  We  get  a  graded   Lie   algebra
structure on $\Di$ but, since now
\bea
\lna  D,ab\rna&=&D(ab)=D(a)\cdot   b+(-1)^{\la a,D+(i-1)\za\ran}a\cdot
D(b)-D(\1)\cdot ab=\\
&&\lna  D,a\rna b+(-1)^{\la a,D+(i-1)\za\ran}a\cdot\lna  D,b\rna-  \lna
D,\1\rna\cdot ab,
\eea
instead of the Leibniz rule we get its generalization.

\begin{theorem} There is a unique $(n+1)$-graded Jacobi bracket
$\lna\cdot, \cdot\rna$ of degree $(-\za,-1)$ on the commutative  algebra
$\Di(\A)$  of  the  previous   theorem, satisfying   (\ref{bra}).   The
associative subalgebra  $\De(\A)$  is  a  graded  Jacobi  subalgebra  of
$\Di(\A)$.
\end{theorem}
The above theorem tells that uniform elements  $D$  from
$\Di_i$  have the degree $(\vert D\vert+(i-1)\za,i-1)$  with  respect  to
the bracket. The  generalized  Leibniz rule reads
\be\label{LR}
\lna A,B\cdot C\rna=\lna A,B\rna\cdot C+(-1)^{\la A+(i-1)\za,B+j\za\ran
+(i-1)j}B\cdot\lna A,C\rna-\lna A,\1\rna\cdot B\cdot C
\ee
and the properties of the graded bracket, written explicitly, are
\bea\label{ss}
\lna  A,B\rna&=&-(-1)^{\la A+(i-1)\za,B+(j-1)\za\ran+(i-1)(j-1)}\lna   B,A
\rna,\\ \label{ji}
\lna\lna A,B\rna,C\rna&=&\lna A,\lna B,C\rna\rna
-(-1)^{\la A+(i-1)\za,B+(j-1)\za\ran+(i-1)(j-1)}\lna B,\lna A,C\rna\rna.
\eea

\medskip\noindent
{\bf Remark.} It is obvious by definitions that if $\za$ and $\za'$ have
the  same  parity,  i.e.   $\za-\za'$   has   even   coordinates,   then
$\Di_i(\A)$ coincides  with  ${\cal D}^{\za'}_i(\A)$  and  the  graded
Jacobi algebras $\Di(\A)$ and ${\Di}'(A)$ are isomorphic.

\medskip
An  element $S\in\Di_2(\A)$ in the     graded     Jacobi      algebra
$(\Di(\A),\cdot,\lna\cdot,\cdot\rna)$  we  call  a  {\it supercanonical
structure} in $\A$ if
\begin{description}
\item{(i)} $\la\vert S\vert+\za,\vert S\vert+\za\ran$ is an even number,
and
\item{(ii)} $\lna S,S\rna=0$.
\end{description}

Similarly as in \cite{Kr2} one proves that any  supercanonical  structure
$S$ in $\A$ determines a graded Jacobi bracket $\{\cdot,\cdot\}_S$ in the
graded algebra $\A$ by
\be\label{gjb}
\{ a,b\}_S=(-1)^{\la a+\za,S+\za\ran}S(a,b).
\ee
We have changed slightly the original definition by a sign in  order  to
get the proper  graded  Lie  algebra  bracket.  This  time  the  bracket
$\{\cdot,\cdot\}_S$ is a Jacobi and  not  Poisson  bracket  in  view  of
(\ref{LR}). Note that if the degree of $S$ is just $\za$, then (i) is
satisfied automatically  and  $\{  a,b\}_S=S(a,b)$.  Thus  we  have  the
following.

\begin{theorem} For  any  supercanonical structure  $S\in\Di(\A)$  (resp.
$S\in\De(\A)$) the formula (\ref{gjb})  defines  a  graded Jacobi (resp.
graded  Poisson)  bracket   on   $\A$   of    degree   $\vert   S\vert$.
Conversely,   every    graded    Jacobi    (resp.    Poisson)    bracket
$\{\cdot,\cdot\}$ of degree $\za$ on $\A$ determines by $S(a,b)=\{ a,b\}$
a supercanonical structure $S$ of $\Di(\A)$  (resp.  $\De(\A)$).
\end{theorem}

We  will  call  this  supercanonical structure   $S\in\Di(\A)$   to   be
{\it associated} with the graded Jacobi bracket of degree $\za$.
Having  a  graded  Lie   algebra   structure   on   $\Di(\A)$   and    a
supercanonical structure  $S\in\De(\A)$ (resp. $S\in\Di(\A)$) we  have,
in  a  standard way,  a cohomology operator $\pa_S=\lna S,\cdot\rna$
which maps $\De_i(\A)$   (resp.   $\Di_i(\A)$)    into    $\De_{i+1}(\A)$
(resp. $\Di_{i+1}(\A)$). This operator has the form of  a  `hamiltonian
vector field'. The corresponding cohomology we will  denote  by
$H^*_P(\A,S)$ (resp. $H^*_J(\A,S)$).
For any Poisson (resp. Jacobi) bracket of degree  $\za$  on  $\A$   we
have  then  the cohomology   operator   $\pa_S$   and   the corresponding
cohomology $H^*_P(\A,S)$ (resp. $H^*_J(\A,J))$, for the associated
supercanonical structure  $S\in\Di(\A)$, which we will call  the  {\it
Poisson} (resp. {\it Jacobi}) {\it cohomology} of  the graded Poisson
(resp. Jacobi) algebra $(\A,\cdot,S)$.

We can deform canonically these cohomology operators in  the  spirit  of
E.~Witten \cite{Wi} as follows.

\begin{lemma} If $S\in\Di(\A)$ is a  supercanonical structure of  degree
$\za$ then
\begin{description}
\item{(i)} $\lna S,S(\1)\rna=0$;
\item{(ii)} $\lna S,S(\1)\cdot A\rna=-S(\1)\cdot A$.
\end{description}
\end{lemma}
\begin{pf} By definition, $0=\lna S,S\rna(\1)=2\lna S,S(\1)\rna$ and, in
view of the generalized Leibniz rule,
\be
\lna S,S(\1)\cdot A\rna=\lna S,S(\1)\rna\cdot A-S(\1)\cdot\lna S,A\rna -
S(\1)\cdot S(\1)\cdot A=S(\1)\cdot\lna S,A\rna,
\ee
due to (i) and $S(\1)\cdot S(\1)=0$.
\end{pf}

\medskip\noindent
It follows from the above lemma that
\be
\pa^t_{LJ}(A)=\lna S,A\rna+tS(\1)\cdot A
\ee
has square 0, i.e. it is a cohomology operator of degree $(\za,1)$ on
$\Di(\A)$ for any parameter  $t$.  The  operator  $\pa^0_{LJ}$  is  just
$\pa_S$ and $\pa^1_{LJ}$, denoted simply $\pa_{LJ}$,  we  will  call  the
{\it Lichnerowicz-Jacobi cohomology operator}. The associated cohomology
will be denoted by $H^*_{LJ}(\A,S)$ and called {\it  Lichnerowicz-Jacobi
cohomology} of the algebra $\A$ associated  with  the  Jacobi  structure
$S$.

\medskip
Note that the general Jacobi  brackets  as  above  do  not
always split into biderivation and derivation, as in the classical case.
There is a new interesting  feature  that  there  exist  graded  Jacobi
brackets being  bidifferential  operators  of order  0.

\medskip\noindent
{\bf Example 3.}  Consider  the  following  bracket   defined   on   the
Grassmann  algebra $\zW(TM)$ of standard differential forms on $M$ by
$\{\za,\zb\}_\zm=\za\we\zm\we\zb$, where $\zm$  is   a   fixed   1-form.
This bracket is clearly of order 0 as a bidifferential operator and of
graded degree 1. Moreover,
\be
\{\za,\zb\}_\zm=-(-1)^{(\za+1)(\zb+1)}\{\zb,\za\}_\zm
\ee
and any double bracket $\{\{\za,\zb\}_\zm,\zg\}_\zm$ vanishes, so that the
Jacobi  identity   is   automatically   satisfied.   In   other   words,
$\{\cdot,\cdot\}_\zm$ is a graded  Jacobi  bracket  on  the  algebra  of
differential  forms  of  order  0   and   degree   1.   In   particular,
$\zm=\{\1,\1\}_\zm$.

\medskip
In fact, every graded Jacobi  bracket $\{\cdot,\cdot\}=S$ of $n$-degree
$\za\in\Z^n$ admits a more general decomposition
\bea\label{zzz}
\{ a,b\}&=&\zL(a,b)+(\zG\cdot I)(a,b)+(c\cdot I^2)(a,b)=\\
&&\zL(a,b)+a\zG(b)-(-1)^{\la a+\za,\za\ran}\zG(a)b+2acb.
\eea
for certain $\zL\in\De_2$, $\zG\in\De_1$, $c\in\De_0=\A$.
Explicitly,
\be
\zL=S-S(\1)\cdot I+\frac{1}{2}S(\1,\1)\cdot I^2,\quad
\zG=S(\1)-S(\1,\1)\cdot I,\quad c=\frac{1}{2}S(\1,\1).
\ee
Here, clearly, $I^2=I\cdot I$.
This can be generalized as follows.

\begin{theorem} Every $D\in\Di_i$  splits into
\be\label{dec}
D=D_0+\frac{1}{1!}D_1\cdot I+\cdots +\frac{1}{i!}D_i\cdot I^{i},
\ee
where
\be\label{dec1}
D_l=\sum_{k=0}^{i-l}(-1)^kD\!\underbrace{\!(\1,\dots,\1)}_{(k+l)-
\textrm{times}}\cdot\frac{I^{k}}{k!}\in\De_l(\A).
\ee
\end{theorem}

\begin{pf} Note that, due to  graded  commutativity,  $I\cdot  I=0$  and
$D(\1,\1)=0$ in
the case when $\la\za,\za\ran$ is even. In the case when $\la\za,\za\ran$
is odd, one shows easily by induction that
\be
(A\cdot I^{n})(\1)=A(\1)\cdot I^{n}+nA\cdot I^{(n-1)}.
\ee
It is now just a matter  of direct calculations  to  show  that
$D_l(\1)=0$ and, using the  identity
\be
\sum_{k=0}^s\frac{1}{k!(s-k)!}=0
\ee
for $s>0$, that the equation (\ref{dec}), with $D_l$'s defined by
(\ref{dec1}), is tautological.
\end{pf}

\medskip
\begin{corollary}.   The   decomposition   (\ref{dec})   determines   an
identification of $\Di_i(\A)$ with
\begin{description}
\item{(i)} $\De_i(\A)\oplus\De_{i-1}(\A)$ in the case  when  $\la\za,\za
\ran$ is even;
\item{(ii)} $\De_i(\A)\oplus\De_{i-1}(\A)\oplus\dots\oplus\De_0(\A)$  in
the case when $\la\za,\za\ran$ is odd.
\end{description}
\end{corollary}
\begin{pf} In the case when $\la\za,\za\ran$ is even, $I\cdot  I=0$,  so
that the splitting (\ref{dec}) reduces to two terms. In the other  case,
all $I^{k}$ are non-zero, due to the formula $I^{k}(a)=
kI^{(k-1)}\cdot a$ which can be easily proved by induction.
\end{pf}

\medskip\noindent
Using these decompositions one can describe  the  bracket  in  $\Di$  in
terms of the bracket in $\De$. In the  case  when  $\la\za,\za,\ran$  is
even it is completely analogous to  the  ungraded  case  (see  the  next
example). For the other case it  is  sufficient  to  use  the  following
lemma.

\begin{lemma} Suppose that  $\la\za,\za,\ran$ is odd. Then,
\be
\lna A\cdot  I^n,B\cdot  I^m\rna=(-1)^{\la(n-1)\za,B\ran}(m(i+1)-n(j+1))
A\cdot B\cdot I^{n+m-1}+(-1)^{\la n\za,B\ran}\lna A,B\rna\cdot I^{n+m}.
\ee
\end{lemma}

\medskip\noindent
{\bf Example 4.} For the (ungraded) algebra $\A=C^\infty(M)$ the  graded
Poisson algebra $\De(\A)={\cal DE}(\A)$ is simply the Gerstenhaber algebra
$A(TM)$  of multivector fields with the (reversed) wedge product  and
the  Schouten bracket. The graded  Jacobi  algebra  $\Di(\A)={\cal
D}(\A)$  is  in  this  case  the Grassmann   algebra
$A(TM\oplus\R)=\textrm{Sec}(\bigwedge(TM\oplus\R))$
with the bracket described in \cite{GM}. Since  here  $I\cdot  I=0$,  we
have the identification
\be
{\cal DE}_i(\A)={\cal DE}_i(\A)\oplus{\cal DE}_{i-1}(\A)=A^i(TM) \oplus
A^{i-1}(TM).
\ee
The Schouten-Jacobi bracket reads (cf. \cite{GM}, formula (27)):
\bea
&&\lna   A_1+I\we   A_2,B_1+I\we   B_2\rna=[A_1,B_1]+(-1)^aI\we[A_1,B_2]+
I\we[A_2,B_1]\\
&&+aA_1\we B_2-(-1)^abA_2\we B_1+(a-b)I\we A_2\we B_2
\eea
for $A\in\Di_{a+1}(\A)$, $B\in\Di_{b+1}(\A)$.
Here $[\cdot,\cdot]$ is the Schouten-Nijenhuis bracket and we use the
standard wedge product. For a Jacobi bracket $S=\zL+I\we\zG$ this  gives
the  cohomology   operator:
\be
\pa_S((B_1,B_2))=([\zL,B_1]+\zL\we B_2+b\zG\we  B_1,-[\zL,B_2]+[\zG,B_1]
+(1-b)\zG\we B_2).
\ee
This is  exactly  (up  to  differences  in  conventions  of   signs  for
the Schouten bracket)  the
cohomology operator of 1-differentiable  Chevalley-Eilenberg  cohomology
of the Jacobi bracket on $C^\infty(M)$,  introduced  and  considered  in
\cite{Li}. The Lichnerowicz-Jacobi cohomology operator is in this case
\be
\pa_{LJ}((B_1,B_2))=([\zL,B_1]+\zL\we B_2+(b+1)\zG\we
B_1,-[\zL,B_2]+[\zG,B_1] -b\zG\we B_2).
\ee
This cohomology has been extensively studied in \cite{LMP,LLMP}  and  we
refer  to \cite{LLMP}  for  more  particular  examples   and   explicit
calculations of particular cohomology.

\medskip\noindent
{\bf Example 5.} Let $\A=A(E)$ be the Grassmann algebra of a vector bundle
$E$. The graded Jacobi bracket (\ref{jb}) on $A(E)$  corresponds  to a
supercanonical structure $S\in{\cal  D}^{-1}(A(E))$.  The  decomposition
(\ref{zzz}) in this case reads
\be
S=S_0+i_\zF\cdot Deg-i_\zF\cdot I,
\ee
where  $S_0$  is  the  supercanonical  structure  corresponding  to  the
Schouten-Nijenhuis bracket.

\medskip\noindent
{\bf Remark.} Suppose we start not from a graded algebra but  from  just
an $n$-graded vector space $V$. We can define inductively  spaces  ${\bf
A}_i^\za(V)$ of multilinear maps in $V$ similarly   to   $\Di_i(\A)$
just relaxing the assumption on the generalized Leibniz  rule.  Then  we
can define a  graded  Lie  bracket $\lna\cdot,\cdot\rna^{NR}$ on  ${\bf
A}^\za(V)=\oplus_{i=0}^\infty{\bf A}^\za_i(V)$ completely along the same
lines as the Jacobi
bracket on $\Di(\A)$.  It  is  just  a  graded  Lie  bracket  satisfying
(\ref{ss}) and (\ref{ji}) but the Leibniz rule has no meaning.  We  will
call this bracket the {\it Nijenhuis-Richardson bracket}  of multilinear
maps in $V$. For a graded algebra $\A$ the space $\Di(\A)$ is just a Lie
subalgebra in ${\bf A}^\za(\A)$. It can be shown that  this  is  exactly
the Nijenhuis-Richardson bracket defined in \cite{LMS} for the graded
vector space $V$ with the grading shifted by $\za$.

\section{Conclusions}
The notion of a Jacobi algebra turns out to unify  various  concepts  of
algebra and differential geometry. A particularly interesting case  is
a linear Jacobi bracket on the exterior algebra of a vector bundle, i.e.
a Jacobi algebroid.
Using  a  gauging  method  we  have  obtained  general
cohomology and homology theories  which  include  a  whole  spectrum  of
(co)homology associated with classical Poisson  and  Jacobi  structures.
Along these lines  we  have  extended  also  the  notion  of  a  Courant
algebroid. Our concept agrees with the known generalizations for the  case
of the tangent bundle.

The construction of a  canonical  Jacobi  algebra
$\Di(\A)$ associated  with  a  given  graded  commutative  algebra  $\A$
generalizes Schouten brackets  and  allows  to  view  the  corresponding
cohomology operators as `hamiltonian vector fields' of a  supercanonical
structure.

One  can  develop  further  this   algebraic    theory
defining    an appropriate
dual object to $\Di(\A)$ (forms)  and  the  corresponding  Cartan-Jacobi
differential calculus. We postpone these problems to a separate paper.

\bigskip\noindent
{\bf Acknowledgments.} The authors address their thanks to Professors
Y. Kosmann-Schwarzbach, I.~S. Krasil'schchik and J.~C. Marrero for
useful comments and suggestions that have helped to shape the final form
of this paper.

\end{document}